\documentclass{article}

\usepackage{xypic}

\usepackage[centertags]{amsmath}
\usepackage{amssymb}
\usepackage{latexsym}
\usepackage{amscd}
\usepackage{theorem}
\usepackage{epsfig}

{\theorembodyfont{\slshape}
        \newtheorem{thm}{Theorem}[section]
        
        \newtheorem{lem}[thm]{Lemma}
        \newtheorem{prop}[thm]{Proposition}
}

{\theorembodyfont{\rmfamily}
        \newtheorem{defn}[thm]{Definition}

        \newtheorem{exa}[thm]{Example}
        
        \newtheorem{notation}[thm]{Notation}
        \newtheorem{assumption}[thm]{Assumption}

}

\makeatletter

\renewcommand{\subsection}{\@startsection{subsection}{3}%
        {\z@}{-3.25ex plus -1ex minus-.2ex}{-1em}{\bf}}
\makeatother

\newcommand{\proof}{\noindent {\bf Proof:\ }}
\newcommand{\Endproof}{\hspace*{\fill} $\Box$ \vspace{1ex} \noindent }

\newcommand{\ZZ}{\mathbb{Z}}

\newcommand{\CC}{\mathbb{C}}
\newcommand{\PP}{\mathbb{P}}

\newcommand{\FF}{\mathbb{F}}

\renewcommand{\O}{{\mathcal O}}

\newcommand{\Spec}{\mathop{\rm Spec}} 
\newcommand{\Hom}{\mathop{\rm Hom}\nolimits}

\newcommand{\Ord}{{\rm ord}}

\newcommand{\dR}{{\rm\scriptscriptstyle dR}}

\newcommand{\Der}{\mathop{\rm Der}\nolimits}
\newcommand{\End}{\mathop{\rm End}\nolimits}
\newcommand{\ord}{{\rm ord}}
\newcommand{\MC}{\mathop{\rm MC}}
\newcommand{\Div}{\mathop{\rm div}}
\newcommand{\Frob}{{\rm Frob}}

\title{Asymptotically good towers and differential equations}
\author{Peter Beelen and Irene I.\ Bouw} \date{}

\begin{document} \maketitle

%19.4 in een file

\begin{abstract}
This paper concerns towers of curves over a finite field with many
rational points, following Garcia--Stichtenoth and Elkies. We present
a new method to produce such towers. A key ingredient is the study of
algebraic solutions to Fuchsian differential equations modulo $p$. We
apply our results to towers of modular curves, and find  new
asymptotically good towers.

 \end{abstract}

\section{Introduction}\label{introsec}
Let $p$ be a prime and $q=p^a$, for some $a>0$. Consider a projective
smooth curve $X$ of genus $g$, defined over $\FF_q$, and write
$N_q(X)$ for the number of $\FF_q$-rational points of $X$. We write
$N_q(g)$ for the maximum of $N_q(X)$, taken over all curves $X$ of
genus $g$ which are defined over $\FF_q$.  The
Drinfel'd--Vl\u{a}du\c{t} bound \cite{DrinfeldVladut} states that
\[
A(q):=\limsup_{g\to\infty}\frac{N_q(g)}{g}\leq \sqrt{q}-1.
\]
Moreover, $A(q)> c\log(q)$, where $c>0$ is a constant
\cite{Serre83}.

 Garcia--Stichtenoth \cite{GS} constructed many examples of infinite
towers of curves $\cdots \to X_{m+1}\to X_{m}\to\cdots \to X_0$
defined over a finite field $\FF_q$ such that the limit of
$N_{q^2}(X_m)/g(X_m)$ is $q-1$. Such towers are called {\sl
asymptotically optimal}. If this limit is positive, the tower is
called {\sl asymptotically good}. Asymptotically optimal towers have, for
example, interesting applications to coding theory
(\cite{Goppa,tsfa-vlad}).  For these applications it is important to
have explicit equations for the curves $X_m$. Garcia--Stichtenoth
define towers of curves recursively, starting from a correspondence
$(g,h):X_0\rightrightarrows X_{-1}$, by taking suitable (normalized)
Cartesian products. A nice feature of this recursive definition is
that one obtains explicit equations for all curves $X_m$ starting from
an equation for $(g,h)$.  One has to choose the correspondence very
carefully for the corresponding tower to have many rational
points. Garcia--Stichtenoth find correspondences that work, but they
do not give a systematic method for finding such correspondences.

Elkies (\cite{Elkies97}) applies this approach to correspondences
$(g,h): X_0(\ell^2)\rightrightarrows X_0(\ell)$, for certain small
values of $\ell$. Here $g$ is the natural projection and $h$ is the
composition of $g$ with the Atkin--Lehner involutions on both
sides. The corresponding tower is $\cdots \to X_0(\ell^{m+1})\to
X_0(\ell^m)\to \cdots $.  This gives equations for modular curves
$X_0(\ell^m)$ starting from equation from
$(g,h):X_0(\ell^2)\rightrightarrows X_0(\ell)$. This is a second
important application of the theory.

Elkies also constructs other asymptotically optimal towers of curves,
starting from correspondences between other Shimura varieties, such as
Drinfel'd modular curves.  These Shimura varieties are moduli spaces
of curves, surfaces etc, and the correspondences are an analog of the
Hecke correspondences for modular curves. Elkies shows that all
asymptotically optimal towers constructed by Garcia--Stichtenoth et.\
al.\ are of this form (\cite{Elkies97}, \cite{Elkies01}). Elkies
suggests that all asymptotically optimal towers arise in this way
(\cite[Fantasia]{Elkies97}).

In this paper we give a new method for constructing asymptotically
good towers. We extract the essential ingredients
from the approach of Garcia--Stichtenoth and Elkies, and formulate a
general set-up. Our approach is concrete and not just applicable to
towers of modular curves. This allows for a more systematic search for
asymptotically optimal towers.

We start from a correspondence $(g,h): X_0\rightrightarrows X_{-1}$
over a finite field $\FF_q$, together with a Fuchsian differential
equation on $X_{-1}$. We say that the correspondence $(g,h)$ is
adapted to the differential equation if the pull back via $g$
is equivalent to the pull back  via $h$ (see Section \ref{desec}
for precise definitions). The correspondence $(g,h)$ gives rise to a
tower of curves ${\mathcal T}_{g,h}=(X_m)_{m\geq 0}$. Under some technical
assumptions we show the following.

\bigskip\noindent {\bf Theorem} \ref{thm:asgood}
 {\sl The tower
${\mathcal T}_{g,h}$ is asymptotically good. This means that the limit
of $N_{q}(X_m)/g(X_m)$ is positive.}

\bigskip\noindent One of the assumptions we make is that $g$ and $h$
are tame, i.e.\ the characteristic of the ground field does not divide
the ramification indices. We also give a criterion for the tower
${\mathcal T}_{g,h}$ to be asymptotically optimal (Theorem
\ref{thm:asopt}).

The reason why such towers have many rational points is roughly the
following.  We suppose that the differential equation has an algebraic
solution $\Phi$. After extending the field of definition $\FF_q$ of
the correspondence, we may assume that the zeros and poles of $\Phi$
are $\FF_q$-rational. The set of these zeros and poles has a subset
${T}$ with the following property. For every $P\in
{T}$ and every $m$, the inverse image of $P$ in $X_m$
consists of unramified and $\FF_q$-rational points.

 It appears that all known examples of tame asymptotically optimal
towers can be reformulated in these terms. The reason is that, by
Elkies' work, the known examples come from certain correspondences
between Shimura curves. (In fact, the tame towers are all towers of
modular curves.) Such moduli spaces come naturally equipped with a
differential equation: the Picard--Fuchs differential equation of a
versal family of the objects it parameterizes.  For towers of modular
curves we work this out in Section \ref{modularsec}. We expect that it
is possible to generalize (parts of) our method to wildly ramified
towers.

The idea for using differential equations for studying the growing
behavior of rational points in a tower came from \cite{GS}. In that
paper Gau\ss' hypergeometric differential equation was used to prove a
property for the Deuring polynomial. We show that the arguments of
\cite{GS} vastly simplify and generalize if one makes a  more
systematic use of differential equations. Our method is also related
in spirit to older work of Ihara (see \cite{Ihara} for a
survey). However, Ihara's work only applies to towers of Shimura
curves. Moreover, it uses $p$-adic uniformization to count points. We
work purely in characteristic $p$ which is more convenient in practice.

  To find new examples of asymptotically good towers, we construct
 correspondences via pull back. Given a correspondence
 $(g,h):X_0\rightrightarrows X_{-1}$ and an arbitrary map $f:Y_{-1}\to
 X_{-1}$, we define a new correspondence $(\tilde{g},
 \tilde{h}):Y_0\rightrightarrows Y_{-1}$. This gives a systematic
 construction of the towers of modular curves found by Elkies in
 \cite[Appendix]{LMS}. This allows to find very many asymptotically
 good towers.

The situation for asymptotically optimal towers is more
complicated. We give a criterion for the pull back of an
asymptotically optimal tower to be again asymptotically optimal
(Theorem \ref{thm:pullbackopt}).  Since our approach does not use the
interpretation of the curves we consider as Shimura varieties, one
might expect to find counter examples to Elkies' Conjecture. However,
we did not find such an example. The reason is that in Theorem
\ref{thm:pullbackopt} there is one condition which is hard to control.
In a later paper we will come back to the question whether this idea
can be used as evidence for Elkies' Conjecture.

The organization of the paper is as follows. In Section \ref{desec} we
review and extend known results on Fuchsian differential equations on
curves in positive characteristic. In Section \ref{pointsec} we give
the recursive definition of a tower of curves corresponding to a
correspondence and establish basic properties. We estimate how the
genus and the number of rational points grow in the tower, and prove
the criterion for a tower to be asymptotically good. In Section
\ref{pullbacksec} we develop the construction of correspondences via
pull back and construct  new examples. Section \ref{modularsec}
reformulates and extends some results of Elkies on towers of modular
curves.

\section{Fuchsian differential equations}\label{desec}

In this section we recall some standard results on Fuchsian
differential equations. For proofs and more details we refer to
\cite[Section 11]{Katz} and \cite{Honda}. Let $k$ be a field of
characteristic $p>0$ and $X/k$ a smooth projective curve. Let $K=k(X)$
be the function field of $X$. Suppose that $M$ is a finite-dimensional
vector space over $K$.

\begin{defn}\label{connectiondef}
A $k$-connection $\nabla$ on $M$ is an additive map
\[
\nabla:M\to \Omega^1_{K/k}\otimes_K M
\]
satisfying the Leibniz rule
\[
\nabla(fm)={\rm d}f\otimes m+f\nabla(m),
\]
for $f\in K$ and $m\in M$.
\end{defn}

Equivalently (\cite[Section 1.0]{Katz}), $\nabla$ corresponds to a
$K$-linear map \[ \nabla:\Der(K/k)\to \End_k(M) \] such that
$\nabla(D)(fm)=D(f)m+f\nabla(D)m$ for $D\in {\Der}(K/k)$, $f\in
  K$ and $m\in M$.  A {\sl horizontal morphism} from $(M_1, \nabla_1)$
  to $(M_2, \nabla_2)$ is a morphism $ \varphi:M_1\to M_2$ of $K$-vector
  spaces which is compatible with the connections, i.e.\
  $\varphi(\nabla_1(D) m)=\nabla_2(D)(\varphi(m))$. We write $\MC(X)$
  for the category of $K$-modules with connection.

\begin{defn}\label{singdef}
Let $P$ be a place of $K/k$ and $t$ a local parameter at $P$. For a
 $K$-basis ${\mathbf e}$ of $M$, we write $\nabla({\rm d}/{\rm d}t)
 {\mathbf e}=A_{\mathbf e} \cdot {\mathbf e}$ with $A_{\mathbf e}\in
 M_n(K)$, where $n=\dim_K M$. We say that $P$ is a {\sl singular point}
 of $(M, \nabla)$ if the matrix $A_{\mathbf e}$ has a pole at $P$ for
 every basis ${\mathbf e}$. 
\end{defn}

\begin{defn}\label{cyclicdef}
We say that $(M, \nabla)$ is {\sl cyclic} if there exists a vector
$m\in M$ and a nonzero derivation $D\in {\Der}(K/k)$ such
that $m, \nabla(D)m, \ldots \nabla^{d-1}(D)m$ span $M$ over $K$, where $d=\dim_K M$.
\end{defn}
It is shown in \cite[11.4]{Katz} that the notion of a cyclic module is
independent of the choice of the derivation $D$. In the rest of this
section, we suppose that $M$ is cyclic, and  of $K$-dimension $d=2$.
A $K$-basis $m, \nabla(D)m$ of $M$ is called a {\sl cyclic basis}.

If $P$ is place of $K/k$, we write ${\mathcal O}_P$ (resp.\
${\mathfrak m}_P$) for the local ring (resp.\ the maximal ideal) at
$P$. Let
\[
\Der_P(K/k)=\{D\in {\Der}(K/k)\, |\,
D({\mathfrak m}_P)\subset {\mathfrak m}_P\}.
\]
If $t=t_P$ is a local parameter of $X$ at $P$, then ${\Der}_P(K/k)$ is
  a free ${\mathcal O}_P$-module with basis $t\,{\rm d/d}t$.

\begin{defn}\label{regsingdef}
Let $P$ be a singular point of $(M,\nabla)$ and $t$ a local parameter
at $P$. For a $K$-basis ${\mathbf e}$ of $M$, we write $\nabla(t\,{\rm
d}/{\rm d}t) {\mathbf e}=B_{\mathbf e} \cdot {\mathbf e}$ with
$B_{\mathbf e}\in M_2(K)$.  We say that $P$ is a {\sl regular
singularity} if there exists a $K$-basis ${\mathbf e}$ of $M$ such
that $B_{\mathbf e}$ is holomorphic at $P$. If all singularities of
$(M,\nabla)$ are regular, we say that $(M, \nabla)$ is a {\sl Fuchsian
module}.

Suppose that $P$ is a regular singularity of $(M,\nabla)$, and let
$B_{\mathbf e}\in M_n({\mathcal O}_P)$ be as above. Write $B_{\mathbf
e}(0)$ for the value of $B_{\mathbf e}$ at $t=0$. The characteristic
polynomial of $B_{\mathbf e}(0)$ is called the {\sl indicial
  equation}. Its roots are the {\sl local exponents.}
\end{defn}

If $P$ is not a singularity its local exponents are $0,1=\dim_K
M-1$. The converse need not be true. Singularities with local
exponents $0,1=\dim_K M-1$ are called {\sl apparent
singularities}.

We now associate to $(M, \nabla)$ a $2$nd order differential equation.
Let $P$ be a regular singularity of $(M, \nabla)$ and $t$ a local
parameter at $P$.  Let
${\mathbf e}=(e_1, e_2:=\nabla({\rm d}/{\rm d}t)(e_1))$ be a cyclic basis
of $M$.  Write
\begin{equation}\label{cyclicbasiseq}
\nabla(\frac{\rm d}{{\rm d}t}){\mathbf e}=A\cdot {\mathbf e}, \qquad\mbox{
with }\qquad A=\begin{pmatrix}0&-a_2\\
1&-a_1\end{pmatrix} .
\end{equation}
It is easy to check that the fact that $P$ is a regular singularity
means that we may choose ${\mathbf e}$ such that $a_i$ has a pole
of order at most $i$ at $P$ for $i=1,2$.

Let $M^\ast=\Hom_{K}(M, K)$ be the $K$-linear dual of $M$. We define a
$k$-connection $\nabla^\ast$ on $M^\ast$ by requiring that
\[
\langle \nabla(D) (m), m^\ast\rangle+\langle m,
\nabla^\ast(D)(m^\ast)\rangle=D(\langle m, m^\ast\rangle),
\]
for $D\in \Der(K/k)$, $m\in M$ and $m^\ast\in M^\ast$. One easily
checks that $\nabla^\ast$ is a connection. In fact, with respect to
the dual basis ${\mathbf e}^\ast$ of $M^\ast$, we have
$\nabla^\ast({\rm d}/{\rm d} t){\mathbf e}^\ast=-A^t{\mathbf
e}^\ast$. Here $A^t$ is the transpose of $A$.

Write $\hat{M}^\ast_P=M^\ast\otimes_K\hat{K}_P$, where $\hat{K}_P$ is the completion of $K$ at $P$.  Denote by
  $(\hat{M}^\ast_P)^{\nabla^\ast}$ the horizontal sections, i.e.\
\begin{equation}\label{deeq}
f_1e_1^\ast+f_2e_2^\ast\in (\hat{M}^\ast_P)^\nabla\qquad \mbox{
  iff }\qquad
\left\{\begin{array}{l}
f_2=f_1',\\
L(f_1):=f_1''+a_1 f_1'+a_2f_1=0.
\end{array}\right.
\end{equation}
This is the {\sl differential equation corresponding to} $(M,
\nabla)$. Giving $(M,\nabla)$ is equivalent to giving the differential
equation $(\ref{deeq})$. We sometimes call $(M, \nabla)$ itself a
differential equation.

One computes that
\[
\nabla(t\frac{{\rm d}}{{\rm d}t})(e_1, te_2)=
\begin{pmatrix}0&-t^2a_2 \\
1& 1-ta_1
\end{pmatrix}(e_1, te_2).
\]
Writing $a_i =c_i t^{-i}+t^{-1}(\cdots)$ with $c_{i}\in k$, we find
that the indicial equation is
\begin{equation}\label{indicialeq}
X^2+(-1+c_1)X+c_2=0.
\end{equation}
 The local exponents $\gamma_1, \gamma_2$ are the roots of this
equation. Note that our notion of local exponents agrees with the
classical ones. This is the reason for taking the differential
equation corresponding to the horizontal sections of $\hat{M}_P^\ast$
rather than $\hat{M}_P$.

Our next topic is algebraic solutions of Fuchsian differential
equations in positive characteristic, following Honda
\cite{Honda}. Let $(M,\nabla)\in \MC(X)$ be a cyclic module of
dimension $2$, and let $P\in X$ be a regular singularity with local
parameter $t$. Choose a cyclic basis ${\mathbf e}$ of $M$. Let $L/K$
be an algebraic extension.  We say that $u\in L$ is an {\sl algebraic
solution} of $(M,\nabla)$ if it is a solution of the corresponding
differential equation (\ref{deeq}). This is equivalent to the fact
that $ue_1^\ast+u'e_2^\ast\in ({ M}^\ast\otimes_K L)^{\nabla^\ast}$. In what
follows we mainly consider solutions in $K$.

\begin{prop}\label{degsolprop}
Let $(M,\nabla)\in \MC(X)$ be a cyclic module of dimension $2$, and let
$u\in K$ be  an algebraic solution (with respect to some choice of a
cyclic basis ${\mathbf e}$).
\begin{itemize}
\item[(i)] Suppose that $P$ is a regular singularity. Write $\gamma_1,
 \gamma_2$ for its local exponents. Then
\[
\ord_P(u)\equiv \gamma_i\bmod{p},
\]
for some  $i\in\{1,2\}$. In particular, $\gamma_i\in \FF_p$.
\item[(ii)] If $P\in X$ is a zero of $u$ we have
\[
\ord_P(u)\equiv 1\bmod{p}\qquad \mbox{ or }\qquad \ord_P(u)\equiv 0\bmod{p}.
\]
\end{itemize}
\end{prop}

\proof Let $u$ be an algebraic solution. Suppose that $P$ is a regular
singularity. Choose a local parameter $t$ at $P$.  Put
$\delta:=\ord_P(u)$. In the complete
local ring $\hat{\O}_P$, we may write $u=t^\delta(\sum u_i t^i)$ with
$u_0$ a unit of $K$. By assumption, $u$ satisfies
\[
 u''+a_1 u'+a_2 u=0,
\]
where $a_i $ has a pole of order at most $i$, since $P$ is a regular
singularity. Write
\[
 a_1=c_{1}t^{-1}+\cdots,\qquad
a_2=c_{2}t^{-2}+\cdots .
\]
Substituting this in the differential equation and taking the
coefficient of $t^{\delta-2}$, we find that
$\delta^2+\delta(c_1-1)+c_2$. Since the indicial equation
(\ref{indicialeq}) is $X^2+(c_1-1)X+c_2$, (a) follows.

If $P$ is a regular point, we may suppose that $a_i $ does not have
a pole at $P$ for $i=1,2$. Hence (b) immediately follows from the
differential equation.  \Endproof

\begin{prop}\label{polysolprop}
Let $(M,\nabla)\in \MC(X)$ be a cyclic module of dimension $2$, and
let $u_1, u_2\in K$ be algebraic solutions (with respect to some
choice of a cyclic basis ${\mathbf e}$). We suppose that
$\ord_P(u_1)\equiv \ord_P(u_2)\bmod{p}$ for some regular singularity
$P$. We write $\Div(u_i)$ for the divisor of $u_i$ on $X$. Then
\[
\Div(u_1)\equiv \Div(u_2)\bmod{p}.
\]

\end{prop}

\proof
This proof follows \cite[Proposition 5.1]{Honda}.

 Consider the set $\Sigma_\delta$ of algebraic solutions $u\in K$ of
$u''+a_1u'+a_2u=0$ which are holomorphic at $P$ and such that
$\ord_P(u)\equiv \ord_P(u_1)\bmod{p}$.

Let $v_1\in \Sigma_\delta$ be a solution whose order at $P$ is minimal
and write $\delta:=\ord_P(v_1)$. Suppose there exists $v_2\in
\Sigma_\delta$ such that $\Div(v_1)\not \equiv \Div(v_2)\bmod{p}$. It
is no restriction to suppose that $v_2$ has minimal order at $P$ among
all algebraic solutions with this property. Write
$\ord_P(v_2)=\delta+p\nu$. Then there exists a constant $c\in \bar{k}$
such that the order of $w:=v_2-ct^{p\nu}v_1\in \Sigma_\delta$ at $P$
is strictly less than $\delta+p\nu$. This contradicts the choice of
$v_2$. Therefore every $u\in \Sigma_\delta$ differs by a $p$th power
from $v_1$. This proves the proposition.
\Endproof

\begin{exa}\label{GMexa}
A key example of a Fuchsian differential equation we will be
interested in in this paper, is the one coming from the Gau\ss--Manin
connection on the modular curve $X(2)$. We recall the situation from
\cite{Katz84}. The statements are easy to generalize to other modular
curves (Section \ref{modularsec}). Let $S=\Spec(\ZZ[\lambda,
1/2\lambda(\lambda-1)])$ and write ${\mathcal E}\to S$ for the
elliptic curve over $S$ given by $y^2=x(x-1)(x-\lambda)$.

We denote by $M:=H^1_\dR({\mathcal E}/S)$ the first de Rham
cohomology group, and by $\nabla:M\to
\Omega^1_S\otimes M$ the Gau\ss--Manin
connection. Write
\[
\omega=\frac{{\rm d} x}{y}=\frac{ {\rm d}
x}{[x(x-1)(x-\lambda)]^{1/2}},\qquad
\omega':=\nabla(\frac{\partial}{\partial \lambda})\, \omega= \frac{{\rm
d} x}{2[x(x-1)]^{1/2}(x-\lambda)^{3/2}}.
\]
Then $\omega$ and $\omega'$ form a basis of $M$. One computes that
\[
\nabla(\frac{\partial}{\partial \lambda})(\omega, \omega')=
\begin{pmatrix}
0&\displaystyle{ -\frac{1}{4\lambda(\lambda-1)}}\\
1&\displaystyle{ -\frac{2\lambda-1}{\lambda(\lambda-1)}}
\end{pmatrix}(\omega, \omega').
\]
The corresponding differential equation (\ref{deeq}) is
\begin{equation}\label{Gausseq}
\lambda(\lambda-1)u''+(2\lambda-1)u'+\frac{1}{4}u=0.
\end{equation}

The differential equation (\ref{Gausseq}) is Gau\ss ' hypergeometric
differential equation. It has three singularities $0,1,\infty$ with
local exponents $0,0;0,0;1/2,1/2$.  Working out the statement of
Proposition \ref{degsolprop} for the singularity $P=\infty$, we obtain
the following. Let $u\in k(\lambda)$ be an algebraic solution. After
multiplying $u$ with a $p$th power, we may suppose that $u$ is a
polynomial. Then $\deg(u)\equiv -\gamma_i\bmod{p}$, were $\gamma_1,
\gamma_2$ are the local exponents at $\infty$. In our case we find
therefore that $\deg(u)\equiv -1/2\bmod{p}$.

The Deuring polynomial (or Hasse
invariant)
\begin{equation}\label{Hasseeq}
\Phi:=\sum_{i=0}^{(p-1)/2}\binom{(p-1)/2}{i}^2\lambda^i
\end{equation}
is a solution $\bmod{p}$ of this differential equation of degree
$(p-1)/2$. Proposition \ref{polysolprop} implies that every other
algebraic solution in characteristic $p$ is of the form $\psi^p \Phi$, for some $\psi\in K$.
\end{exa}

In general, a module $(M,\nabla)\in \MC(X)$ does not have algebraic
solutions. Honda \cite[appendix]{Honda} shows that $(M,\nabla)$ has
``sufficiently many solutions in a weak sense'' if and only if the
$p$-curvature of $(M,\nabla)$ is nilpotent. The notion of sufficiently
many solutions in a weak sense is stronger than just the existence of
an algebraic solution; we refer to \cite{Honda} for a
definition. However, if $(M,\nabla)$ is cyclic, has dimension two, and
its singularities are $\FF_p$-rational than the two notions are
equivalent \cite[Cor.\ 1 to Prop.\ 2.3]{Honda}. Katz \cite{Katz} shows
that in a ``geometric'' context (like in Example \ref{GMexa}) the
$p$-curvature is always nilpotent, in particular, the corresponding
differential equation has an algebraic solution in some extension
 $L/K.$

Suppose that $(M_1,\nabla_1)$ and $(M_2, \nabla_2)$ are elements of
$\MC(X)$. We define a $k$-connection $\nabla$ on $M=M_1\otimes_K M_2$ by
putting $\nabla(D)(m_1\otimes m_2)=(\nabla_1(D)m_1)\otimes
m_2+m_1\otimes (\nabla(D)m_2)$, for $D\in \Der(K/k)$ and $m_i\in M_i$.

\begin{defn}\label{equivdef}
Let $(M_1, \nabla_1)$, $(M_2, \nabla_2)\in \MC(X)$ be modules with
  $\dim_K M_1=\dim_K M_2$. We say that $(M_1, \nabla_1)$ is {\sl
  equivalent} to $(M_2, \nabla_2)$ is there exists a one-dimensional
  module $(M_3, \nabla_3)\in \MC(X)$ such that:
\begin{itemize}
\item $(M_1, \nabla_1)\otimes_K(M_3, \nabla_3)\simeq (M_2, \nabla_2)$,
\item $(M_3, \nabla_3)$ has an algebraic solution $\theta$,
\item the set of singularities of $(M_3, \nabla_3)$ is contained in
  the set of singularities of $(M_1, \nabla_1)$.
\end{itemize}
\end{defn}

In terms of local coordinates this definition means the following. Let
$(M_1, \nabla_1)\in \MC(X)$ be a cyclic, Fuchsian module of dimension
two and let $(M_3, \nabla_3)\in \MC(X)$ be a one-dimensional Fuchsian
module. Let $P$ be a regular singularity of both $(M_1, \nabla_1)$ and
$(M_3, \nabla_3)$ with local parameter $t$. Choose a cyclic basis
$(e_1, e_2)$ for $M_1$ as in (\ref{cyclicbasiseq}), i.e.\ we write
\[
\nabla_1({\rm d}/{\rm d}t)(e_1, e_2)=
\begin{pmatrix}0&-a_2 \\1&-a_1 \end{pmatrix}(e_1, e_2).
\]

We identify $M_3$ with $K$, and write $\nabla_3({\rm d}/{\rm
d}t)=B\cdot 1. $ Then with respect to the basis $\xi_1=e_1\otimes 1$
and $\xi_2=Be_1\otimes 1+e_2\otimes 1$ of $M_2=M_1\otimes_K M_3$ we
find
\[
\nabla_2({\rm d}/{\rm d}t)(\xi_1, \xi_2)=
\begin{pmatrix} 0&\displaystyle{-a_2
    +B a_1+B'-B^2}\\
1&\displaystyle{-a_1+2B}
\end{pmatrix}(\xi_1, \xi_2).
\]
The corresponding differential equation is
\[
y''+(a_1-2B)y'+(-B'-B a_1+a_2+B^2)y=0.
\]
One checks that if $u$ is an (algebraic) solution of the differential
equation corresponding to $(M_1, \nabla_1)$ then $\theta u$ is an
(algebraic) solution of $(M_2, \nabla_2)$. Here $\theta$ is an algebraic solution of $M_3$ (Definition \ref{equivdef}).
Let $\gamma$ (resp\ $\gamma_1, \gamma_2$) be the local exponents of
$(M_3, \nabla_3)$ (resp.\ $(M_1, \nabla_1)$) at $P$. One computes that
the indicial equation of $(M_2,\nabla_2)$ at $P$ is
$(X-\gamma_1-\gamma) (X-\gamma_2 -\gamma)$.

Suppose we are given $(M,\nabla)\in \MC(X)$ and a cover $f:Y\to X$
defined over $k$, i.e.\ $f$ is a finite separable map between smooth
and absolutely irreducible curves. Let $L=k(Y)$ be the function field
of $Y$.

\begin{defn}\label{pullbackdef}
We define the {\sl pull back} $( M_f, \nabla_f)$ on $Y$ as
follows. Write $\nabla(m)=m\otimes{\rm d}g$, where ${\rm d} g\in
\Omega^1_K$. Then $ M_f=M\otimes_K L$ and $ \nabla_f:M\otimes_K L\to
\Omega^1_L\otimes (M\otimes_K L)$ is defined by $m\otimes 1\mapsto
f^\ast(m\otimes {\rm d}g):=m\otimes {\rm d} (g\circ f)\in M\otimes_K
L\otimes_L\Omega^1_L=M_f\otimes_L\Omega^1_L.$
\end{defn}

In local coordinates this may be described as
follows. Let $Q$ be a point of $Y$ and $P$ its image in $X$. Choose a
local parameter $s$ of $Q$ and let $t=s^e$ be a local parameter of
$P$. Here $e$ is the ramification index of $Q$ in $f$. Write $f'(s)\in
\O_Q$ for the derivative of $f$ at $Q$. Choose an appropriate basis
${\mathbf e}=(e_1, e_2)$ of $M$ at $P$, and write
\[
\nabla({\rm d}/{\rm d}t){\mathbf e}=
\begin{pmatrix}0&-a_2 \\ 1&-a_1 \end{pmatrix}{\mathbf e},
\]
as above. Then
\[
\nabla_f({\rm d}/{\rm d}s)(e_1, f'(s)e_2)=
\begin{pmatrix}0&\displaystyle{-(f'(s))^2 a_2(f(s))}\\
1&\displaystyle{-f'(s)a_1(f(s))+f''(s)/f'(s)}
\end{pmatrix}(e_1, f'(s)e_2).
\]
The corresponding differential equation is
\[
L_f(v)=\left(\frac{\rm d}{f'(s){\rm d}s}\right)^2
v+a_1(f(s))\left(\frac{\rm d}{f'(s){\rm
    d}s}\right)v+a_2(f(s)) v=0.
\]

One easily checks that if $(M,\nabla)$ is Fuchsian, then $(M_f,
\nabla_f)$ is Fuchsian also. Moreover, with notation as above, if $P$
is a regular singularity with local exponents $(\gamma_1, \gamma_2)$
then $(e\gamma_1, e\gamma_2)$ are the local exponents at $Q$. Note
that it may happen that $P$ is a singularity but $Q$ is not. (It is
easy to characterize this in terms of the local monodromy, but we do
not need this here.)

\begin{notation}\label{polysolnot}
Let $(M, \nabla)\in \MC(X)$ and write $S$ for its set of
singularities. Suppose that $(M, \nabla)$ has an algebraic solution
$\Phi$.  Let $f:Y\to X$ be a cover and $(M_f, \nabla_f)$ the pull back
module. We write $\Phi_f:=\Phi\circ f$ for the corresponding algebraic
solution of $(M_f,
\nabla_f)$.
\end{notation}

\begin{defn}\label{corrdef}
Let $(M, \nabla)\in \MC(X)$.  A {\sl
correspondence adapted to} $(M, \nabla)$ is a pair of (separable)
covers $g, h:Y\rightrightarrows X$ between smooth and absolutely
irreducible curves such that the pull back modules $(M_g, \nabla_g)$
and $(M_h, \nabla_h)$ on $Y$ are equivalent. The correspondence is
{\sl trivial} if there exists an automorphism
$\sigma:X\stackrel{\sim}{\to}X$ such that $h=\sigma\circ g$. A
correspondence $(g,h)$ is called {\sl tame} if the covers $g$ and $h$
are tame.
\end{defn}

A correspondence $(g,h):Y\rightrightarrows X$ may equivalently be
described by giving a curve $C\subset X\times X$. Here $C=\{( g(P),
h(P))\, |\, P\in Y \}$ is the {\sl curve of correspondence.} The {\sl
degree of the correspondence } is the cardinality of $\{ (x, y)\in
C\}$, where $x\in X$ is a fixed, sufficiently general point. We will
be particularly interested in correspondences of degree one. In this
case the map $Y\to C$ defined by $P\mapsto (g(P), h(P))$ is
generically an isomorphism.

\begin{prop}\label{feprop}
Let $(M,\nabla)\in \MC(X)$ be a Fuchsian, cyclic module of dimension
two. Let $S$ be its set of singularities. Suppose that
\begin{itemize}
\item $(M,\nabla)$ has an algebraic solution $\Phi$, and
\item there exists a point $P\in S$ whose local exponents $\gamma_1,
\gamma_2$ are equal.
\end{itemize}
 Let $(g,h):Y\rightrightarrows X$ be a correspondence adapted to
$(M,\nabla)$. Write ${\mathfrak S}$ for the set of singularities of
the pull back module $(M_g, \nabla_g)$. Then there exists a divisor
$D=\sum_{P_i\in {\mathfrak S}} n_i P_i$ such that
\[
\Div(\Phi_g)\equiv \Div(\Phi_h)+D\bmod{p}.
\]
\end{prop}

\proof Note that ${\mathfrak S}$ is also the set of singularities of
$(M_h, \nabla_h)$. The fact that $(M_g, \nabla_g)$ and $(M_h,
\nabla_h)$ are equivalent means that there exists a one-dimensional
module $(N, \nabla_N)$ such that $(M_g, \nabla_g)\otimes (N,
\nabla_N)\simeq (M_h, \nabla_h)$. It is no restriction to suppose that
we have equality. Recall that there exists an algebraic function
$\theta$ such that $\theta\Phi_g$ is a solution of $(M_h,
\nabla_h)$. Moreover, the poles and zeros of $\theta$ are in
${\mathfrak S }$. This function is an algebraic solution of the module
$(N, \nabla_N)$. Put $D:=\Div(\theta)$.

Let $P$ be as in the statement of the proposition and let $Q$ be a
point of $Y$ with $g(Q)=P$. Write $e_g$ (resp.\ $e_h$) for the
ramification index of $Q$ in $g$ (resp.\ $h$).  Let $\gamma$ be the local
exponent of $(N, \nabla_N)$ at $P$. By comparing the local exponents,
we see that the sets $\{e_g\gamma_1+\gamma, e_g\gamma_2+\gamma\}$ and 
$\{e_h\gamma_1, e_h\gamma_2\}$ are equal. By assumption, $\gamma_1=
\gamma_2$. Therefore Proposition \ref{degsolprop} implies that
$\ord_Q(\theta\Phi_g)\equiv e_g\gamma_1+\gamma\equiv e_h\gamma_1\equiv
\ord_Q(\Phi_h)\bmod{p}$. The statement now follows from Proposition 
\ref{polysolprop}.
\Endproof

\begin{lem}\label{pullbacklem}
Let $(M, \nabla)\in \MC(X)$ and $g,h:Y\rightrightarrows X$ be a
correspondence adapted to $(M, \nabla)$. Then, for every map
$\phi:Z\to Y$, the modules $(M_{g\circ\phi}, \nabla_{g\circ\phi})$ and
$(M_{h\circ\phi}, \nabla_{h\circ\phi})$ are equivalent also.
\end{lem}

\proof Straight forward.
\Endproof

\section{Estimates for the number of points and the genus in a tower} \label{pointsec}

In this section we define a tower of curves from a tame correspondence
adapted to a differential equation $(M, \nabla)$. We also estimate the
genus (Proposition \ref{genusprop}) and number of points (Proposition
\ref{countprop}) in the tower. The results are easiest to understand
in the well-known case of towers of modular curves (Section
\ref{modularsec}). It may be helpful to look at this case before
reading the proofs in the general case.

Let $(M,\nabla)$ be a Fuchsian differential equation of rank $2$ with
set of singularities $S$. We always suppose that $M$ is cyclic
(Definition \ref{cyclicdef}). We denote by $(g,h): X_0
\rightrightarrows X_{-1}$ a tame correspondence adapted to
$(M,\nabla)$ (Definition \ref{corrdef}) unbranched outside $S$. We
always assume that $X_0$ and $X_{-1}$ are smooth and absolutely
irreducible curves. We assume that the covers $g$ and $h$ are disjoint
(i.e.\ the covers $g,h:X_0 \to X_{-1}$ do not have a common subcover,
or alternatively, there do not exist functions $\phi, \psi_1$ and
$\psi_2$ such that $g=\phi\circ \psi_1$ and $h=\phi\circ\psi_2$ and
$\deg \phi \neq 1$). Denote the common set of singularities of
$(M_g,\nabla_g)$ and $(M_h,\nabla_h)$ by $\mathfrak S$. To the
correspondence $(g,h)$ we associate a tower of curves
\begin{equation}\label{towereq} {\mathcal T}_{g,h}=\left(
\cdots\stackrel{\displaystyle\pi_m}{\longrightarrow}X_m\stackrel{\displaystyle
\pi_{m-1}}{\longrightarrow} X_{m-1}\to \cdots
\stackrel{\displaystyle\pi_{0}}{\longrightarrow}X_0\right),
\end{equation} where $X_m$ is a smooth projective curve and
$\pi_m$ is a cover. For $m \ge 1$, the curve $X_m$ is the
normalization of the curve $X_m'$. The curves $X_m'$ are
defined recursively as follows. We put $X_0'=X_0$. For $m\geq 0$,
we define $X_{m+1}'$ by the fiber product
\[
\xymatrix{
X_{m+1}'\ar[d]_{\pi_{m}}\ar[r]&X_0\ar[d]_h\\
X_{m}'\ar[r]_{g_{m}}&X_{-1},}
\]
where $\pi_m(x_0,\ldots,x_{m+1})=(x_0,\ldots,x_{m})$ and $g_m(x_0,
\ldots, x_{m})=g(x_{m})$. We have

$$X_m'=\{(x_0,\dots,x_m) \in X_0^{m+1} \, | \, h(x_i)=g(x_{i-1}),
\, 1 \le i \le m\}.$$ The cover $\pi_m: X_{m+1}' \to X_m'$ induces
a cover from $X_{m+1}$ to $X_m$ which we again denote by $\pi_m$.
We will also denote by $x_i: X_m \to X_0$ (with $0 \le i \le m$)
the cover defined by $x_i: P \mapsto x_i(P)$.

To obtain a tower of curves, we require that $X_m$ is an absolutely
irreducible curve, for all $m$. Clearly, a necessary condition for
this is that $g$ and $h$ are disjoint. This condition is not
sufficient however.  For example if we take $g(t)=t^2+t$ and
$h(t)=1/(t^2+t)$ both defined over a finite field of characteristic
two, then $g$ and $h$ are disjoint, but the corresponding curve $X_2$
is not irreducible (although $X_1$ is).

We now state a sufficient condition for all curves $X_m$ occurring in
the tower ${\mathcal T}_{g,h}$ to be absolutely irreducible.  Let $Y$
and $X$ be curves defined over a field $k$ and suppose we are given a
cover $\pi: Y \to X$ defined over $k$. We say a point $P$ of $X$ is
{\sl totally branched} in the cover $\pi$ if there exists a point $Q$
of $Y$ with $\pi(Q)=P$ such that $e(Q|P)=\deg \pi$. The following
lemma is obvious.

\begin{lem}\label{lem:irred} Suppose that the cover $\pi_m: X_m
\to X_{m-1}$ has a totally branched point for any $m \ge 1$. Then all
curves $X_m$ are absolutely irreducible. \end{lem}

The condition that $g$ and $h$ are disjoint is not vital for the
construction of the tower. If $g=\phi\circ
 \psi_1$ and $h=\phi\circ\psi_2$, one should replace the
 correspondence $(g,h)$ by $(\psi_1, \psi_2)$. Note that if $g$ and $h$ are
 disjoint we have $\deg \pi_m=\deg h$. All asymptotically
 good towers of function fields found by Garcia--Stichtenoth et al.\
 (see e.g.\ \cite{GS1,GS,GS2,GS3}) can be described as a tower $\mathcal
 T_{g,h}$ for a suitable correspondence $(g,h)$. For the definition of the tower, we do not need that the
correspondence is adapted to some differential equation. We only
need this afterwards to estimate the number of $\FF_q$-rational
points.

We now state some restrictions we will assume in the rest of this
section.

\begin{assumption}\label{notation}

\begin{itemize}

\item[(a)]{$\mathfrak S=g^{-1}(S)=h^{-1}(S)$.}

\item[(b)]{All curves $X_m$ occurring in the tower ${\mathcal T}_{g,h}$
are absolutely irreducible.}

\item[(c)]{$\deg g=\deg h=:\delta$.}

\end{itemize}

\end{assumption}

We will usually check (b) by using Lemma \ref{lem:irred}. Note
that (c) is a natural restriction, since if $\deg g \neq \deg h$
the tower ${\mathcal T}_{g,h}$ is asymptotically bad \cite{skew}.

We start by estimating the genus $g(\mathcal T_{g,h})$ of a tower
$\mathcal T_{g,h}$. This genus is defined in the following way:

$$g(\mathcal T_{g,h}):=\lim_{m \to \infty} \frac{g(X_m)}{\delta^m}.$$
Here $g(X_m)$ denotes the genus of the curve $X_m$. This limit
exists \cite{GS3}, but may be infinite. A necessary condition for
a tower $\mathcal T$ to be asymptotically good is that $g(\mathcal
T)<\infty$. The following proposition checks this in our
situation.

\begin{prop}\label{genusprop} Let $(M,\nabla)$ be a Fuchsian
differential equation of rank $2$ with set of singularities $S$.
Suppose Assumption \ref{notation} holds. Then for any $m$ and any
point $P$ of $X_m$ we have 
\begin{itemize}
\item[i)]{$x_{m-1}(P)\in\mathfrak S \Longleftrightarrow
x_{m}(P)\in\mathfrak S,$} 
\item[ii)]{$g({\mathcal T}_{g,h}) \le
g(X_0)+\frac{\# {\mathfrak S}-2}{2}.$} 
\end{itemize}

\end{prop}

\proof We extend the constant field to ${\overline{\FF}_q}$,
which does not make a difference since we are only interested in the
genus at this point. We show that the branch locus of the tower
${\mathcal T}_{g,h}$ is contained in $\mathfrak S$. (By branch locus we
mean here the set of points of $X_0$ that are branched in the cover $X_m
\to X_0$, for some $m$.)

Let $P$ be a point of the curve $X_m$. By the recursive definition
of the tower, we have $h(x_m(P))=g(x_{m-1}(P))$. Therefore, if
$x_{m-1}(P) \in \mathfrak S$, we have $x_m(P) \in
h^{-1}g(x_{m-1}(P)) \subset h^{-1}g(\mathfrak S)=\mathfrak S$ by
Assumption (a). Conversely, $x_m(P) \in \mathfrak S$ implies
$x_{m-1}(P) \in \mathfrak S$. This proves the first part of the
proposition.

Recall that the following diagram commutes
\[
\xymatrix{
X_{m}\ar[d]_{\pi_{m-1}}\ar[r]&X_0\ar[d]_h\\
X_{m-1}\ar[r]_{g_{m-1}}&X_{-1}.} 
\] 
If $P \in X_m$ ramifies in $\pi_{m-1}:X_m \to X_{m-1}$ then
$g(x_{m-1}(P)) \in S$, since we assumed that $h$ is unbranched outside
$S$. We distinguish two cases: $x_m(P) \in \mathfrak S$ and $x_m(P)
\not\in \mathfrak S$. If $x_m(P) \in \mathfrak S$, one obtains from
the first part by induction $x_0(P) \in \mathfrak S$. Now assume
$x_m(P) \not\in \mathfrak S$. Let $e$ be the ramification index of
$x_m(P)$ in the cover $h:X_0 \to X_{-1}$. Since $x_m(P)$ is a regular
point of $(M_h,\nabla_h)$, its local exponents are $0,1$. By
part i) and the assumption $x_m(P) \not\in \mathfrak S$ we have
$x_{m-1}(P) \not\in \mathfrak S$. By considering the local exponents,
we conclude that $x_{m-1}(P)$ has ramification index $e$ in the cover
$g:X_0 \to X_{-1}$. By Abhyankar's lemma, we conclude that the cover
$\pi_0: X_{1} \to X_{0}$ is unbranched at $x_{m-1}(P)$.  Consider the
commutative diagram
\[
\xymatrix{
X_{m}\ar[d]_{\pi_{m-1}}\ar[r]^{\psi}&X_1\ar[d]_{\pi_0}\\
X_{m-1}\ar[r]_{x_{m-1}}&X_0,}
\]
with $\psi:X_m \to X_1$ induced by the map $\psi': X_m' \to X_1'$
defined by $\psi(P)=(x_{m-1}(P),x_m(P))$ and $x_{m-1}: X_{m-1} \to
X_0$ defined by $x_{m-1}:Q \mapsto x_{m-1}(Q)$. It follows that
the cover $\pi_{m-1}$ is unramified at $P$, i.e., $x_0(P)$ does
not belong to the branch locus of the tower.

One uses the Riemann--Hurwitz formula for the cover $X_m \to X_0$ to deduce
\[
\frac{g(X_m)-1}{\delta^m} \le g(X_0)-1+\frac{\#\mathfrak S \cdot
(\delta^m-1)}{2\delta^m}.
\]
%\left(\frac{1}{\delta^m}+\cdots+\frac{1}{\delta}\right).$$
%The right-hand side of this equation evaluates to
%$g(X_0)+(\delta^m-1)(\# \mathfrak S-2)/2 \delta^m$.
The proposition follows by letting $m$ tend to infinity. \Endproof

We now investigate the asymptotic behavior of the number of
rational points in the tower $\mathcal T_{g,h}$. A key role is is
played by Proposition \ref{feprop}. Let $\Phi$ be an algebraic
solution of $(M,\nabla)$. Further, let $\Phi_g$ (resp.\ $\Phi_h$)
denote the corresponding solution of $(M_g,\nabla_g)$ (resp.\
$(M_h,\nabla_h)$) (Notation \ref{polysolnot}).

Let $\pi:Y\to X$ be a cover of curves over $k$ and $P$ is a
$k$-rational point of $X$. We say that $P$ is {\sl completely split}
if $P$ is unbranched and every point $Q$ of $Y$ with $\pi(Q)=P$ is
$k$-rational. The following set will turn out to describe a set of
completely split places of $X_0$ in the tower $\mathcal
T_{g,h}$. Define
\begin{equation}\label{Teq}
\mathfrak T:=\{x_0 \in X_0 \, | \, \ord_{x_0}\Phi_g \not\equiv 0
\bmod{p} \makebox{ and } x_0 \not\in \mathfrak S\}.
\end{equation}
Recall that Proposition \ref{degsolprop} implies that $\Ord_{x_0}\Phi_g\equiv 1\bmod{p}$ for $x_0\in \mathfrak{T}$. 
 The following lemma gives some properties of this set. It will be
useful in the investigation of the number of rational points in the
tower.

\begin{lem}\label{felem} Let $(M,\nabla)$ be a Fuchsian
differential equation of rank $2$ with set of singularities $S$.
Suppose Assumption \ref{notation} holds. Further let $\alpha,\beta
\in X_0$ be such that $h(\beta)=g(\alpha)$. Then
\begin{itemize}
\item[i)]{$\alpha \in \mathfrak T \Longleftrightarrow \beta \in
\mathfrak T,$}
\item[ii)]{$\mathfrak T=\{x_0 \in X_0 \, | \, \ord_{x_0}\Phi_h \not\equiv 0
\bmod{p} \makebox{ and } x_0 \not\in \mathfrak S\}.$}
\end{itemize}
\end{lem}

\proof Suppose $\beta \in \mathfrak T$. By the definition of
$\mathfrak T$ and part i) of Proposition \ref{genusprop}, we have
$\alpha \not\in \mathfrak S$. By Proposition \ref{feprop} we have
$\ord_{\beta} \Phi_h \not\equiv 0 \bmod{p}$. Since
$g(\alpha)=h(\beta)$ and $\alpha\not\in\mathfrak{S}$, we conclude that
$\Phi_h(\beta)=\Phi_g(\alpha)=0$. Proposition \ref{degsolprop} implies
that $\ord_{\alpha} \Phi_g \not\equiv 0 \bmod{p}$. This proves i).

The second part of the lemma follows directly from the (proof of) the
first part.\Endproof

The first part of the above lemma implies that $g(\mathfrak
T)=h(\mathfrak T)$. We write
$$T:=g(\mathfrak T)=h(\mathfrak T).$$
 The second part of Lemma \ref{felem} shows that the role of
$g$ and $h$  in the definition of $\mathfrak
T$ can be interchanged.

Given an absolutely irreducible curve $C$ defined over
$\mathbb{F}_q$, we denote by $N_q(C)$ the number of
$\mathbb{F}_q$-rational points of $C$. For a tower $\mathcal
T_{g,h}=(X_0,X_1,\dots)$ defined as above with constant field
$\mathbb{F}_q$ we define the {\sl splitting rate} of the tower
$\mathcal T_{g,h}$ by
$$\nu_q(\mathcal T_{g,h}):=\lim_{m\to \infty }\frac{N_q(X_m)}{\delta^m}.$$
This limit exists \cite{GS3} and is a nonnegative finite number. A
necessary condition for a tower $\mathcal T$ to be asymptotically
good is $\nu_q(\mathcal T)>0$. The following proposition gives an
estimate for $\nu_q(\cal T)$ in our situation.

\begin{prop}\label{countprop} Let $(M,\nabla)$ be a Fuchsian
differential equation of rank $2$ with set of singularities $S$.
Suppose Assumption \ref{notation} holds. Suppose that the constant
field $\mathbb{F}_q$ of the tower $\mathcal T_{g,h}$ is such that
all points of $X_0$ in the set ${\mathfrak T}$ are
defined over $\mathbb{F}_q$. Then
$$\nu_q(\mathcal T_{g,h}) \ge \# \mathfrak T.$$
\end{prop}

\proof Recall that $\delta:=\deg g=\deg h$. Since $ S$
and $T$ are disjoint, for any $\alpha \in T$ there are exactly
$\delta$ points of $X_0$ lying above $\alpha$.  Moreover, all these
points of $X_0$ are defined over $\mathbb{F}_q$ by our assumption. Let
$P$ be a point of $X_0$ with $g(P)=\alpha$. Write $h(P)=\beta$. Since
$T=g({\mathfrak T})=h({\mathfrak T})$ it follows that $\beta\in T$.

Now suppose that we have constructed inductively $\delta^{m-1}
\#\mathfrak T$ points of $X_{m-1}$ defined over $\mathbb{F}_q$ and
lying above $\mathfrak T$. Consider the commutative diagram
\[
\xymatrix{
X_{m}\ar[d]_{\pi_{m-1}}\ar[r]^{\psi}&X_1\ar[d]_{x_1}\\
X_{m-1}\ar[r]_{x_{m-1}}&X_0,} 
\] 
with $\psi:X_m \to X_1$ induced by
the map $\psi': X_m' \to X_1'$ defined by
$\psi'(P)=(x_{m-1}(P),x_m(P))$, the map $x_1:X_1 \to X_0$ is given by
$Q \mapsto x_1(Q)$, and similarly $x_{m-1}: X_{m-1} \to X_0$ is
defined by $R \mapsto x_{m-1}(R)$. Given an $\alpha \in \mathfrak T$,
we can construct $\delta$ points $P$ of $X_1$ defined over
$\mathbb{F}_q$ with $x_1(P)=\alpha$ and $\delta^{m-1}$ points $Q$ of
$X_{m-1}$ also defined over $\mathbb{F}_q$ with
$x_{m-1}(Q)=\alpha$. Given such a $P$ and $Q$, there exists at least
one point $R$ of $X_{m}$ lying above both $P$ and $Q$. Moreover, by
Lemma \ref{felem}, we have $x_m(R) \in \mathfrak T$. It follows that we
have obtained in this way all $\delta^m$ points of $X_m$ lying above
$\alpha$ and that any of these points is defined over
$\mathbb{F}_q$. \Endproof

The field $\mathbb{F}_q$ mentioned in the above proposition is
called the minimal splitting field of the tower $\mathcal
T_{g,h}$. In other words, we have the following definition.

\begin{defn} Given a correspondence $(g,h):X_0 \to X_{-1}$
defining a tower $\mathcal T_{g,h}$, we define the {\sl minimal
splitting field} of this tower to be the smallest field $k$ such
that \begin{itemize} \item[i)]{the correspondence $(g,h):X_0
\rightrightarrows X_{-1}$ is defined over $k$,} \item[ii)]{all
points of $X_1$ in the set $\pi_0^{-1}(\mathfrak T)$ are defined
over $k$.} \end{itemize} \end{defn}

The following theorem gives a sufficient condition for $\mathcal
T_{g,h}$ to be asymptotically good.

\begin{thm}\label{thm:asgood}
Suppose Assumption \ref{notation} holds. Let $\mathbb{F}_q$ be the
minimal splitting field of $\mathcal T_{g,h}$. Suppose that
$\ord_P(\Phi) \not \equiv 0 \bmod{p}$ for some $P \in \PP^1$ and
that $g^{-1}(P) \not \subset \mathfrak S$. Then the tower
$\mathcal T_{g,h}$ is asymptotically good.
\end{thm}

\proof By Proposition \ref{genusprop}, the tower has finite genus.
We will show that the set $\mathfrak T$ is nonempty. Let $Q \in
g^{-1}(P)$ and $Q \not \in \mathfrak S$. Since all ramification in
the cover $g:X_0\to X_{-1}$ is tame, we have $\ord_Q(\Phi_g) \not
\equiv 0 \bmod{p}$. We conclude that $Q \in \mathfrak T$. By
Proposition \ref{countprop} the tower has positive splitting rate.
Hence the tower $\mathcal T_{g,h}$ is asymptotically good.
\Endproof

\begin{thm}\label{thm:asopt}
Suppose Assumption \ref{notation} holds. Let $\mathbb{F}_q$ be the
minimal splitting field of $\mathcal T_{g,h}$. Suppose that 
\[
2\#
\mathfrak{T}=(\sqrt{q}-1)(\# \mathfrak S+2g(X_0)-2).
\]
 Then the tower $\mathcal T_{g,h}$ is asymptotically optimal.
\end{thm}
\proof This follows immediately from Propositions \ref{genusprop}
and \ref{countprop}. \Endproof

The minimal splitting field is in practice often difficult to
calculate. This is a serious problem in finding asymptotically optimal
towers via the criterion of Theorem \ref{thm:asopt}. Proposition
\ref{zieveprop} is a useful tool to deal with this problem: it
essentially controls the minimal splitting field at the cost of
introducing a new condition on the correspondence
$(g,h):Y\rightrightarrows X$. Namely, we need to suppose that the
correspondence has degree one. In Section
\ref{pullbacksec} we will always make this assumption. In the case of
modular curves (Section \ref{modularsec}) this condition is always
satisfied, see the proof of Lemma \ref{modularzievelem}.

Recall from Section \ref{desec} that if a correspondence
$g,h):Y\rightrightarrows X$ has degree one, then the map $Y\to C$ of
$Y$ onto the curve of correspondence is generically a bijection.

\begin{prop}\label{zieveprop} Let $X$ and $Y$ be smooth and
absolutely irreducible curves defined over $k$, and let
$(g,h):Y \rightrightarrows X$ be a tame correspondence of degree one over $k$.

Let $V \subset X$ be a set of $k$-rational points such that for
any $\alpha, \beta \in Y$ with $h(\beta)=g(\alpha)$ we have
$g(\alpha) \in V \Leftrightarrow g(\beta) \in V$. Let $\alpha \in
Y$ be such that
\begin{itemize}
\item[i)]{$g(\alpha) \in V$.}
\item[ii)]{$g(\alpha)$ is a $k$-rational point of $X$.}
\item[iii)]{$(g(\alpha),h(\alpha))$ is not a singularity of $C$.}
\end{itemize}
Then $\alpha$ is a $k$-rational point of $Y$.
\end{prop}

\proof We first show that $h(\alpha)$ is a $k$-rational point of
$X$. There exists a point $\beta$ such that $h(\alpha)=g(\beta)$.
By the definition of $V$ and $i)$, the point
$h(\alpha)=g(\beta)$ is in $V$ and hence $k$-rational.

Since the map $\phi$ has degree one, it can be inverted for
nonsingular points of $C$. The $k$-rationality of
$(g(\alpha),h(\alpha))$ then implies the $k$-rationality of
$\alpha=\phi^{-1}(g(\alpha),h(\alpha))$. \Endproof

Condition iii) in the above lemma is in practice not a heavy
restriction. Since the number of singularities of $C$ is finite,
they can usually be dealt with by hand in any particular case.
For certain correspondences of degree two this lemma is due to Zieve
(\cite{GS}). We will apply this lemma in the
situation that $V=T(=g(\mathfrak T)=h(\mathfrak T))$. If the
conditions of the above lemma are satisfied, then the points in
the set $\mathfrak T$ are defined over $k$ if the points in the
set $T$ are.

\section{Constructing towers via pull back}\label{pullbacksec}

As before let $(g,h):X_0 \rightrightarrows X_{-1}$ be a correspondence
adapted to a Fuchsian differential equation $(M,\nabla)$, where we
suppose that $g$ and $h$ are disjoint. As always, we suppose that
$(M,\nabla)$ is cyclic. We write $S$ for the set of singularities of
$(M,\nabla)$.  Let $f:Y_{-1} \to X_{-1}$ be a (separable) cover of
smooth, absolutely irreducible curves, which is allowed to  have wild
ramification and may be ramified outside $S$. We suppose that
all curves and covers are defined over a finite field $k$. Write
$(M_f,\nabla_f)$ for the pull back of $(M,\nabla)$ via $f$ and $S_f$
for the set of singularities of $(M_f,\nabla_f)$. In this section we
make the following additional assumption.

\begin{itemize}

\item[(d)]The correspondence $(g,h):X_0\rightrightarrows X_{-1}$ has
degree one.

\end{itemize}

Recall that we defined the curve of correspondence $C\subset
X_{-1}\times X_{-1}$ by
\begin{equation}\label{Cdefeq}
C:=\{(g(P),h(P)) \, | \, P \makebox{ a point of } X_0\}.
\end{equation}
 The curve $C$ is the image of $X_0$ under the map $g*h:X_0 \to X_{-1}
\times X_{-1}$ defined by $(g*h)(P)=(g(P),h(P))$. Assumption (d)
implies that the map $X_0\to C$ has degree one. Denote by $p_1$
(resp. $p_2$) the projections of $C$ onto it first (resp. second)
coordinate. We have the following commutative diagram
\[
\xymatrix{ &&&X_0\ar[d]^{g*h}\ar[dll]_{g}\ar[drr]^{h}\\
&X_{-1}&&C\ar[ll]^{p_1}\ar[rr]_{p_2}&&X_{-1}.}
\]

Let $D\subset Y_{-1}\times Y_{-1}$ be an absolutely irreducible
component of the inverse image of $C$ under $(f,f):Y_{-1} \times
Y_{-1} \to X_{-1} \times X_{-1}$. After extending the field of
definition $k$, we may suppose that $D$ is defined over $k$.  We have
the following diagram
\[
\xymatrix{
&X_0\ar[d]^{g*h}&&&D\ar[dlll]_{(f,f)}\ar[dl]^{\tilde{p}_1}\ar[dr]_{\tilde{p}_2}\\
&C\ar[dl]^{p_1}\ar[dr]_>>>>{p_2}&&Y_{-1}\ar[dlll]_<<<<<<<<<{f}&&Y_{-1}\ar[dlll]_f\\
X_{-1}&&X_{-1}.} \] The maps $\tilde{p}_1$ (resp. $\tilde{p}_2$) are
projections of $D$ onto its first (resp. second) coordinate.  Recall
that we always suppose that the curve $X_0$ is smooth. Denote by $Y_0$
the normalization of the curve $D$, then we have the following diagram
\[
\xymatrix{ &&&Y_0\ar[d]^{}\ar[dll]_{\tilde{g}}\ar[drr]^{\tilde{h}}\\
&Y_{-1}&&D\ar[ll]^{\tilde{p}_1}\ar[rr]_{\tilde{p}_2}&&Y_{-1}.}
\]
The maps $\tilde{g}$ and $\tilde{h}$ are defined such that the
diagram commutes.
\begin{defn}\label{pullback}
We call $(\tilde{g}, \tilde{h})$ the {\sl pull back} of $(g,h)$
under $f$.
\end{defn}

The following lemma gives a key property of the pull back
correspondence $(\tilde{g},\tilde{h})$.

\begin{lem}\label{adaptlem} The correspondence $(\tilde{g}, \tilde{h})$ is adapted
to $(M_{f},\nabla_f)$. \end{lem}

\proof We have the following (commutative) diagram
\[
\xymatrix{ &&&&Y_0\ar[dlll]\ar[dl]^{\tilde{g}}\ar[dr]_{\tilde{h}}\\
&X_0\ar[dl]^{g}\ar[dr]_>>>>{h}&&Y_{-1}\ar[dlll]_<<<<<<<<<{f}&&Y_{-1}\ar[dlll]_{f}\\
X_{-1}&&X_{-1}.}
\]
Hence the result follows immediately from Lemma \ref{pullbacklem}.
\Endproof

Note that if $\deg g=\deg h$, then $\deg \tilde{g}=\deg
\tilde{h}$. The above lemma motivates that if the tower $\mathcal
T_{g,h}$ is asymptotically good, the tower $\mathcal
T_{\tilde{g},\tilde{h}}$ is a good candidate for being asymptotically
good as well. Recall that we defined (\ref{Teq}) a set $\mathfrak T
\subset X_0$ consisting of completely splitting places of the tower
$\mathcal T_{g,h}$. We denote by $T \subset X_{-1}$ the set $g(\mathfrak
T)=h(\mathfrak T)$.

\begin{thm}\label{thm:pullbackgood} Let $(M,\nabla)$ be a Fuchsian
differential equation of rank $2$ and let $(g,h)$ be a correspondence
adapted to $(M,\nabla)$ all defined over a finite field
$\mathbb{F}_q$. Suppose Assumptions (a), (b), (c), (d) hold. Let $f:
Y_{-1} \to X_{-1}$ be a cover and suppose that Assumption (b)
 holds for the pull back correspondence
$(\tilde{g},\tilde{h})$ as well. If the set $\mathfrak T$ is non
empty, then $\mathcal T_{\tilde{g},\tilde{h}}$ is asymptotically good
over some extension field of $\FF_q$. \end{thm}

\proof Our assumptions imply that the ramification locus of the pull
back tower ${\mathcal T}_{\tilde{g}, \tilde{h}}$ is contained in
$f^{-1}(S)$. Therefore Proposition \ref{genusprop} implies that the
genus $g({\mathcal T}_{\tilde{g}, \tilde{h}})$ of the pull back tower
is finite.

We claim that $\nu_q(\mathcal T_{\tilde{g},\tilde{h}})>0$. Denote by
$\phi: Y_0 \to X_0$ the map induced by $f$. Since $\mathfrak T \subset
X_0$ consists of completely splitting points of the tower $T_{g,h}$,
the non empty set $f^{-1}(\mathfrak T) \subset Y_0$ consists of
completely splitting places of the tower $T_{\tilde{g},\tilde{h}}$ if
we extend the constant field suitably. Here we use that ${\mathfrak
T}$ is unbranched in $f$. \Endproof

Before proceeding, we give an example illustrating Theorem
\ref{thm:pullbackgood}.

\begin{exa}\label{goodexa} We consider the correspondence $(g,h):\PP^1
\rightrightarrows \PP^1$ given by $h(t)=t^2$ and $g(t)=4t/(t+1)^2$.
 The correspondence $(g,h)$ is adapted to the Gau\ss' hypergeometric
 differential operator $L(u)=t(t-1)u''+(2t-1)u'+u/4$, which has the
 Deuring polynomial as a solution (Example \ref{GMexa}). We denote the
 Deuring polynomial by $\Phi$. Recall that $S=\{0,1,\infty\}$ is the
 set of singularities of $L$. The set of singularities of the pull
 back differential equation $L_g$ equals ${\mathcal S}=\{0, \pm 1,
 \infty\}$. Therefore the correspondence satisfies Assumptions (a) and
 (c). One checks that the point $x_0=0$ on $X_0$ is totally branched
 in the tower. Therefore Assumption (b) follows from Lemma
 \ref{lem:irred}.

The corresponding tower
${\mathcal T}_{g,h}$ is the tower of modular curves $X_0(2^m)$
starting from $m=3$. The tower $\mathcal{T}_{g,h}$ is essentially the
same as a tower considered in \cite{GS}.

The curve of correspondence
 $C=\{(g(x),h(x)) \, | \, x \in \PP^1\}$ is given by the equation
\begin{equation}\label{Ceq}
4a(b-2)^2-(a+1)^2b^2=0. 
\end{equation}
One checks that $X_0$ is a normalization of $C$. In other words, Assumption (d)
is satisfied. This means that we can apply the pull back construction
to the tower $\mathcal T_{g,h}$. Further note that the point $(-1,2)$
of $C$ is a singularity.

Let $f:\PP^1 \to \PP^1$ be the cover defined by
$t=f(s):=-n(n/(n-1))^{n-1}(s^n-s^{n-1})$ for an integer $n \ge 2$
satisfying $p \not|(n-1)$ and $p\not|n$, where $p$ denotes the
characteristic. In particular $Y_{-1}=\PP^1$. An explicit calculation
shows that the cover $f$ is unbranched outside the set
$\{0,1,\infty\}$. Let $\FF_q$ be the smallest finite field containing
all roots of the polynomial $\Phi(T^m-T^{m-1})$. Using Proposition
\ref{zieveprop} one checks that the pull back tower $\mathcal
T_{\tilde{g},\tilde{h}}$ is asymptotically good over the field
$\FF_q$, for all $n$ for which Assumption (b) holds. We will not
determine the field $\FF_q$ explicitly here. For $n=2$ we obtain an
asymptotically optimal tower which turns out to be the modular tower
$((X_0(2^m))$ starting from $m=4$. For $n>2$ one does not seem to
obtain asymptotically optimal towers. Trivially one sees that $q$
divides $2\cdot n!$.
\end{exa}

Example \ref{goodexa} illustrates how to find asymptotically good
towers via pull back. The only problem is to check Assumption (b) for
the pull back tower. This condition is satisfied if
$\deg(\tilde{g})=\deg(g)$. To obtain asymptotically optimal towers, we
need to impose a condition on the minimal splitting field. This
condition is in practice hard to check.

\begin{thm}\label{thm:pullbackopt} Let $(g,h):X_0\rightrightarrows
X_{-1}$ be a correspondence adapted to a differential equation
$(M,\nabla)$ with singularity set $S$, satisfying Assumptions (a), (b),
(c), (d).  Let $f:Y_{-1} \to X_{-1}$ be a tame cover unbranched
outside $S$. Denote by $(\tilde{g},\tilde{h})$ the pull back
correspondence and suppose it satisfies Assumption (b).
%We demand that $\deg g =\deg \tilde{g}$.

If ${\mathcal T}_{g,h}$ is asymptotically optimal and the towers $\mathcal
T_{g,h}$ and $\mathcal T_{\tilde{g},\tilde{h}}$ have the same minimal
splitting field $\FF_q$, then $\mathcal T_{\tilde{g},\tilde{h}}$
is asymptotically optimal also. \end{thm}

\proof By our assumptions $\phi^{-1}(\mathfrak T)$ consists of
completely splitting places of the tower $\mathcal
T_{\tilde{g},\tilde{h}}$.  Therefore $\nu_q(\mathcal
T_{\tilde{g},\tilde{h}}) \ge \deg \phi \cdot \# \mathfrak T$ by
Proposition \ref{countprop}.

As usual, denote by $\mathfrak S$ the singularities of the
differential equation $(M_g,\nabla_g)$.  The
Riemann-Hurwitz genus formula for the cover $\phi:Y_{0} \to X_{0}$ implies
$$2g(Y_0)-2+\#\phi^{-1}({\mathfrak S})=\deg
\phi\cdot(2g(X_0)-2+\#\mathfrak S).$$ Since $\phi^{-1}({\mathfrak S})$
contains the set of singularities of $(M_{\tilde{g}\circ f},
\nabla_{\tilde{g}\circ f})$, it follows from  Proposition \ref{genusprop}  that
$g(\mathcal T_{\tilde{g},\tilde{h}}) \le \deg
\phi\cdot(g(X_0)+(\#\mathfrak{S}-2)/2)$. Therefore $\lambda(\mathcal
T_{\tilde{g},\tilde{h}}) \ge \lambda(\mathcal T_{g,h})=q-1$ and we are
done.  \Endproof

As a consequence of Theorem \ref{thm:pullbackopt}, we give an
asymptotically optimal tower. We consider again the correspondence
$(g,h):\PP^1 \rightrightarrows \PP^1$ given by
 $h(t)=t^2$ and
$g(t)=4t/(t+1)^2$.
Using Proposition \ref{zieveprop} we immediately obtain that the roots
of the Deuring polynomial are squares in $\FF_{p^2}$. In fact
these roots are fourth powers in $\FF_{p^2}$ (see \cite{GS}).

Let $f:\PP^1 \to \PP^1$ be defined by \[
t=f(s)=\frac{16s^2}{(s-1)^4}. \] One checks that
$S_f=f^{-1}\{0,1,\infty\}=\{0,\pm 1, 3\pm\sqrt{2}, \infty\}$. The pull
back differential equation is \[
L_f(v)=v''+\frac{s^4-4s^3+20s^2+8s-1}{s(s^2-1)(s^2-6s+1)}v'+
\frac{16(s^2+1)}{s(s+1)(s-1)^2(s^2-6s+1)}v=0. \] The pull back of $C$
with respect to the map $(f,f)$ has two absolutely irreducible
components of genus $0$ and one of genus $2$. We write $f(A)=a$ and
$f(B)=b$ and use (\ref{Ceq}). The components of genus $0$ of the pull
back of $C$ are then given by the equations
\begin{equation}\label{genus0} -A+A^2+4AB+B^2-AB^2=0
\qquad\makebox{ and }\qquad
1-A+4AB-AB^2+A^2B^2=0. \end{equation} 
%The component of genus $2$
%satisfies the equation
%
%\begin{equation}\label{genus2} G-2G^2+G^3-4GH+8G^2H-
%4G^3H+H^2-6GH^2+26G^2H^2\end{equation}$$-6G^3H^2+G^4H^2-4GH^3+
%8G^2H^3-4G^3H^3+GH^4- 2G^2H^4+G^3H^4=0.$$

\noindent One may choose any genus $0$ component $D$ from
(\ref{genus0}) and a coordinate $y$ of its normalization $Y_0$
such that the maps $\phi:Y_0\to X_0$ and $\tilde{g},
\tilde{h}:Y_0\to \PP^1$ are described as follows.
\[
\xymatrix{
&X_0\ar[ld]_g\ar[rd]^h&&&Y_0\ar[lll]_\phi\ar[ld]_{\tilde{g}}\ar[rd]^{\tilde{h}}\\
\PP^1&&\PP^1&\PP^1&&\PP^1,} \] with $\phi(y)=4y^2/(y^2-1)^2$,
$\tilde{h}(y)=y^2$ and $\tilde{g}(y)=-y(y-1)/(y+1)$. 

Let ${\mathcal T}_{\tilde{g}, \tilde{h}}=(Y_0,Y_1,\dots,Y_m,\dots)$ be
the tower of curves defined by the correspondence $(\tilde{g},
\tilde{h})$. One checks that $y_0=\infty$ is totally branched in the
tower. Therefore Lemma \ref{lem:irred} implies that the curves $Y_m$
are irreducible for all $m$. Then $Y_{m}$ is given by the equations
\[
y_i^2=-\frac{y_{i-1}(y_{i-1}-1)}{y_{i-1}+1}, \makebox{ with $1 \le
i \le m$}.\]

We write $\Phi_f(y)=\Phi(f(y))$ for the algebraic solution of
$L_f(v)=0$ (Notation \ref{polysolnot}). Using that the correspondence
$(\tilde{g},\tilde{h})$ is adapted to $L_f$ (Lemma \ref{adaptlem}),
one checks that
\[
(y^2-1)^{2p-2}\Phi_f(y^2)=(y^2+1)^{2p-2}\Phi_f\left(\frac{-y(y-1)}{y+1}\right).
\]
This illustrates Proposition \ref{feprop}.

\begin{prop}\label{exaprop} The tower ${\mathcal T}_{\tilde{g},
\tilde{h}}$  is asymptotically optimal if $p\equiv
\pm 1\bmod{8}$. \end{prop}

\proof To apply Theorem \ref{thm:pullbackopt}, we only need to
 determine the minimal splitting field of the tower
$\mathcal
T_{\tilde{g},\tilde{h}}$. Using
Proposition \ref{zieveprop} we see that this field is in fact the
splitting field of $\Phi_f(t)$. In other words, we are interested
in the solutions of the equation
\begin{equation} \label{sol}
\frac{16y^2}{(y-1)^4}=\left(\frac{4y}{(y-1)^2}\right)^2=\lambda \
{\makebox{ with }\ \Phi(\lambda)=0.} \end{equation}
 We have
already seen that all roots of the Deuring polynomial are squares in
$\FF_{p^2}$. Write $\lambda=\mu^2$. Equation (\ref{sol}) has solutions in
$\FF_{p^2}$ if
and only if $\mu+1$ is a square in $\FF_{p^2}$.

Suppose that  $p \equiv \pm 1
\bmod{8}$. We claim that for any root $\lambda$ of $\Phi$ and any
element $\mu$ with $\mu^2=\lambda$ the element $\mu+1$ is a square
in $\FF_{p^2}$. We will prove this claim following the approach by
R{\"u}ck in the appendix of \cite{GS}.

Consider the
elliptic curve $E_{\lambda}$ given by $ Y^2=X(X-1)(X-\lambda).$
Since $\lambda$ is a root of the Deuring polynomial, $E_\lambda$ is
supersingular.
We first suppose that $\lambda \not\in \{-1,2,1/2\}$ and that
$\lambda$ is not a sixth root of unity. It is known that
$\Frob_{p^2}$, the Frobenius automorphism over $\FF_{p^2}$, acts on $E_\lambda$
as multiplication by $\pm p$. This implies that the $x$-coordinate
of any $8$-torsion point of $E_\lambda$ is an element of $\FF_{p^2}$.
Here we use that
$p \equiv \pm 1 \bmod{8}$.

The point $(0,0)$ of $E_{\lambda}$ is a point of order two. For
any point $(a,b)$ satisfying $2(a,b)=(0,0)$ we have
$$a^2=\lambda$$ as can be seen directly from the addition formulas
of $E_{\lambda}$. Hence we may choose $a=-\mu$. For the points
$(c,d)$ satisfying $2(c,d)=(a,b)$ we find in a similar way
$(c-a)^4+4c^2(a-1)^2a=0$. Write $\mu=\nu^2$. Note
that $\nu \in \FF_{p^2}$, since all roots of $\Phi$ are fourth
powers in $\FF_{p^2}$. We obtain
\begin{equation}\label{basiceq}
(c^2+2(-\nu+\mu-\nu \mu)c+\lambda)(c^2+2(\nu+\mu+\nu
\mu)c+\lambda)=0.
\end{equation}
The discriminant of any of these factors is $\mu+1$ up to
multiplication with squares in $\FF_{p^2}$. Since all solutions of
(\ref{basiceq}) are in $\FF_{p^2}$, the claim follows.

If $\lambda \in \{-1,2,1/2\}$, then a direct computation shows
that $\mu+1$ is a square in $\FF_{p^2}$ in our situation. On the
other hand if $\lambda$ is a sixth root of unity, then
$\Frob_{p^6}$, the Frobenius automorphism on $\FF_{p^6}$, acts as
multiplication by $\pm p$ on $E_{\lambda}$. By a similar argument
as above, we conclude that $\sqrt{\mu+1} \in \FF_{p^6}$. On the
other hand, it is obvious that $\sqrt{\mu+1}\in \FF_{p^8}$.
Therefore, $\mu+1$ is a square in $\FF_{p^2}$ in this case as
well.

Theorem \ref{thm:pullbackopt} now implies that the tower $\mathcal
T_{\tilde{g},\tilde{h}}$ is asymptotically optimal. \Endproof

%Using the component given in equation (\ref{genus2}) we can also
%construct a pull back tower $\mathcal T=(Y_0,Y_1,\dots)$ defined
%over the field $\FF_{p^2}$ for any odd prime number $p$. We denote
%by $F$ the polynomial occurring in equation (\ref{genus2}). One
%checks that the function field of the curve $Y_m$ can be described
%by $$\FF_{p^2}(Y_m)=\FF_{p^2}(x_0,x_1,\dots,x_{m+1}) \ \makebox{
%with } \ F(x_{i+1},x_{i})=0 \ \makebox{ for } \  0 \le i \le m.$$
%Analyzing the proof of Proposition \ref{exaprop} we find that for
%this tower $\lambda(\mathcal T) \ge (p-1)/2$ if $p \equiv \pm 1
%\bmod{8}$.

\section{Towers of modular curves}\label{modularsec}
In this section we apply the results of Sections \ref{pointsec} and
\ref{pullbacksec} to towers of modular curves.

Fix an integer $\ell>3$. We do not suppose that $\ell$ is prime.
Write $X_0(\ell^m)$ for the modular curve parameterizing (generalized)
elliptic curves $E$ together with a cyclic isogeny $E\to E'$ of degree
$\ell^m$. For a precise description of the points above $j=\infty$
(the cusps) in terms of generalized elliptic curves we refer to
\cite{DelRap}. The curve $X_0(\ell^m)$ has a natural smooth model
over $\ZZ[1/\ell]$. Denote by $\sigma_m:X_0(\ell^m)\to X_0(\ell^m)$
the Atkin--Lehner involution. It sends an isogeny $E\to E'$ to its
dual isogeny.

We define a correspondence $(g,h):X_0(\ell^2)\rightrightarrows
X_0(\ell)$ as follows. Suppose that $(E_1\to E_2\to E_3)$ corresponds
to point of $X_0(\ell^2)$, i.e.\ $E_1\to E_3$ is a cyclic isogeny of
degree $\ell^2$ and $E_i\to E_{i+1}$ has degree $\ell$. Then $g(E_1\to
E_2\to E_3)=(E_1\to E_2)$ is the standard projection and $h(E_1\to
E_2\to E_3)=(E_2\to E_3)$ is $\sigma_1\circ g\circ \sigma_2$.

Analogous to Example \ref{GMexa}, we obtain a differential equation on
$X_0(\ell^m)$. Fix a prime $p$ relatively prime to $\ell$. We denote
by $X_0(\ell^m)/\FF_p$ the reduction of $X_0(\ell^m)$ to characteristic $p$.

 Let $S=\Spec(\FF_p[j, 1/j(j-1728)])$. Write ${\mathcal E}_{\ell^m}$
for the universal elliptic curve on $X_0(\ell^m)$; it exists for
$\ell^m>3$. Let $M_{\ell^m}=H^1_{\dR}({\mathcal E}_{\ell^m}/S)$ be the
de Rham cohomology group and $\nabla:M_{\ell^m}\to \Omega^1_S\otimes
M_{\ell^m}$ the Gau\ss--Manin connection (\cite[Section
7]{Katz}). Then $(M_{\ell^m},\nabla)\in
\MC(X_0(\ell^m))$; its set of singularities $S_{\ell^m}$ is contained in the
inverse image of $j=0,1728, \infty$.

\begin{lem}\label{modularcorrlem}
The correspondence $(g,h):X_0(\ell^2)\rightrightarrows X_0(\ell)$ is
adapted to $(M_\ell, \nabla)$.
\end{lem}

\proof
 The pull back of ${\mathcal E}_\ell$ via $g$ is just the
universal elliptic curve ${\mathcal E}_{\ell^2}$. Denote the pull back
of ${\mathcal E}_\ell$ via $h$ by ${\mathcal E}_{\ell^2}'$. The
concrete description of $g$ and $h$ given above implies that there is
an isogeny ${\mathcal E}_{\ell^2}\to {\mathcal E}_{\ell^2}'$. It
induces an isomorphism $H^1_{\dR}({\mathcal E}_{\ell^m}/S)\simeq
H^1_{\dR}({\mathcal E}'_{\ell^m}/S)$ on the de Rham cohomology groups.
This implies the statement of the lemma.
\Endproof

On $X(1)$ we still have a versal family of elliptic curves ${\mathcal
  E}_1$. We may choose ${\mathcal E}_1$ such that it is nonsingular
  outside $j=0,1728, \infty$. The differential equation corresponding
  to $(M_1, \nabla)$ is a hypergeometric
  differential equation (i.e.\ a Fuchsian differential equation with
  three singularities). Its singularities are $j=0, 1728, \infty$ with
  local exponents $ 0, 1/3; 0,1/2; 1/12, 1/12$, respectively. Note
  that $(M_{\ell^m}, \nabla)$ is the pull back
  of $(M_{1}, \nabla)$ via the natural
  projection $(E\to E')\mapsto E$. Denote by $\nu_2(\ell^m)$ (resp.\
  $\nu_3(\ell^m)$, resp.\ $\nu_\infty(\ell^m)$) the number of
  singularities of $(M_{\ell^m}, \nabla)$
  above $j=1728$ (resp.\ $j=0$, resp.\ $j=\infty$). Let $\mu(\ell^m)$
  be the degree of $X_0(\ell^m)\to X(1)$. These numbers are computed
  in \cite[Proposition 1.43]{Shimura}. Moreover, it is shown in
  \cite[Proposition 1.40]{Shimura} that
\begin{equation}\label{genusmodulareq}
g(X_0(\ell^m))=1+\frac{\mu(\ell^m)}{12}-\frac{\nu_2(\ell^m)}{4}-
\frac{\nu_3(\ell^m)}{3}- \frac{\nu_\infty(\ell^m)}{2}.
\end{equation}

The {\sl supersingular polynomial} is defined as
\[
\Phi_1(j)=\prod (j-j(E))\in \FF_p[j],
\]
 where the product is taken over the supersingular elliptic curves
 $E/\bar{\FF}_p$. Put
\[
\alpha:=\left[\frac{p}{12}\right],\quad
\delta:=\left\{\begin{array}{ll}
0&\mbox{ if }p\equiv 1\bmod{3}\\
1&\mbox{ if }p\equiv 2\bmod{3}
\end{array}\right.
,\quad
\epsilon:=\left\{\begin{array}{ll}
0&\mbox{ if }p\equiv 1\bmod{4}\\
1&\mbox{ if }p\equiv 3\bmod{4} \end{array}.\right. \] There exists
a polynomial $\tilde{\Phi}_1$ of degree $\alpha$ such that
$\Phi_1=j^\delta(j-1728)^\epsilon\tilde{\Phi}_1$. All zeros of
$\tilde{\Phi}_1$ are simple.

\begin{lem}
The polynomial $\Phi_1$ is an algebraic solution of
$(M_{1}, \nabla)$.
\end{lem}

\proof
This is well known. It can for example be checked by direct
verification, or deduced from \cite{Katz84}.
\Endproof

We denote by $\Phi_{\ell^m}$ the induced algebraic solution of
$(M_{\ell^m}, \nabla)\otimes \FF_p$ (Notation 
\ref{polysolnot}).

\begin{lem}\label{modularzievelem}
We write
\[
{\mathfrak T}_{\ell}:=\{ x\in X_0(\ell)_{\FF_p}\, |\, \Phi_\ell(x)=0,
  \mbox{ and } x\not\in S_\ell\}.
\]
The points of ${\mathfrak T}_{\ell}$ are $\FF_{p^2}$-rational.
\end{lem}

\proof 
It is well-known that the roots of $\Phi$ are rational over
$\FF_{p^2}$ (\cite{Silverman}).

For $j_1,j_2\in X(1)$ we write $j_1\sim_\ell j_2$ if there
exists a cyclic isogeny of degree $\ell$ from the elliptic curve with
$j$-invariant $j_1$ to the elliptic curve with $j$-invariant $j_2$.
Define a (singular) curve
\[
{\mathcal C}(\ell)=\{ (j_1,j_2)\in X(1)\times X(1)\, |\, j_1\sim_\ell
j_2\}.
\]
We obtain a commutative diagram
\[
\xymatrix{
&X_0(\ell)\ar[ddl]_{g_0}\ar[ddr]^{h_0}\ar[d]&\\
&{\mathcal C}(\ell)\ar[dl]^{{\rm pr}_1}\ar[dr]_{{\rm pr}_2}&\\
X(1)&&X(1),
}
\]
with $g_0(E_1\to E_2)=j(E_1)$ and $h_0(E_1\to E_2)=j(E_2)$ (compare to
(\ref{Cdefeq})).

If $E_1$ and $E_2$ are elliptic curves without complex multiplication,
there exists at most one isogeny $E_1\to E_2$ of fixed degree
$\ell$. Namely, $\Hom(E_1, E_2)$ is a right $\End(E_2)$-module of
rank one and $\End(E_2)\simeq \ZZ$, since $E_2$ does not have
CM. This implies that the map $X_0(\ell)\to {\mathcal C}(\ell)$ has
degree one and is defined over $\FF_p$. The lemma now follows from
Proposition
\ref{zieveprop}, since the roots of $\Phi$ are $\FF_{p^2}$-rational.  \Endproof

\begin{prop}\label{modulartowerprop} Let $\ell>3$ be an integer.
\begin{itemize}
\item[(a)] Let $(g,h):X_0(\ell^2)\rightrightarrows X_0(\ell)$ be the
  correspondence defined above. Then the corresponding tower of curves
  is isomorphic to ${\mathcal T}_{g,h}=(X_0(\ell^m))$.
\item[(b)] The tower ${\mathcal T}_{g,h}$ is asymptotically optimal.
\end{itemize}
\end{prop}

\proof Part (a) is proved in \cite{Elkies97}.  Part (b) follows from
the work of \cite{Ihara}. It is also proved in \cite{TVZ}. We indicate an
alternative proof using our results.

 If $\nu_2(\ell^2)=\nu_3(\ell^2)=0$, the proposition follows from
 Theorem \ref{thm:asopt}. Otherwise, the estimates for
 $g(X_0(\ell^m))$ and $N_{p^2}(X_0(\ell^m))$ given in Section
 \ref{pointsec} are not quite good enough. But it is easy to compute
 these quantities directly, using the results of
 \cite{Shimura}. Namely, one checks that
\[
\lim_{m\to\infty} \frac{\nu_2(\ell^{m+1})}{\delta^m}=\lim_{m\to\infty}
\frac{\nu_3(\ell^{m+1})}{\delta^m}=\lim_{m\to\infty}
\frac{\nu_\infty(\ell^{m+1})}{\delta^m}=0.
\]
Therefore (\ref{genusmodulareq}) implies that the genus of the tower is
\begin{equation}\label{genuseq}
\lim_{m\to\infty}
\frac{g(X_0(\ell^{m+1}))}{\delta^m}=\frac{\mu(\ell)}{12}.
\end{equation}

To estimate the splitting rate in the tower, one needs to count the
points on $X_0(\ell^m)$ above $j=0, 1728$ which are not
singularities. Such points above $j=0$ (resp.\ $j=1728$) are zeros of
the pull back of $\Phi$ if and only if $j=0$ (resp.\ $j=1728$) is
supersingular, i.e.\ $p\equiv 2\bmod{3}$ (resp.\ $p\equiv
3\bmod{4}$). Distinguishing cases acording to the value of
$p\bmod{12}$, one finds that
\begin{equation}\label{splittingeq}
\lim_{m\to \infty} \frac{N_{p^2}(X_0(\ell^{m+1}))}{\delta^m}\geq
\frac{(p-1)\mu(\ell)}{12}.
\end{equation}
Equations (\ref{genuseq}) and (\ref{splittingeq}) imply that the tower
is optimal.
\Endproof

It is easy to see that for every $a\geq 1$ we can consider the tower
defined by the correspondence $(g,h):X_0(\ell^{a+1})\rightrightarrows
X_0(\ell^a)$. This yields the subtower $(X_0(\ell^{a+m}))$ which is of
course again asymptotically optimal.

We now present a variant of this construction. Choose an integer
$\lambda$ relatively prime to $\ell$ and $p$. Consider the pull back
of the correspondence $(g,h):X_0(\ell^{a+1})\rightrightarrows
X_0(\ell^a)$ via the natural projection $X_0(\lambda\ell^a)\to
X_0(\ell^a)$. It is easy to see that the pull back correspondence is
$(\tilde{g}, \tilde{h}):X_0(\lambda\ell^{a+1})\rightrightarrows
X_0(\lambda\ell^{a})$.

\begin{prop}
 The tower defined by $(\tilde{g},
\tilde{h}):X_0(\lambda\ell^{a+1})\rightrightarrows
X_0(\lambda\ell^{a})$ is
\[
\cdots\to X_0(\lambda\ell^m)\to \cdots X_0(\lambda\ell^{a+1})\to
X_0(\lambda\ell^{a}).
\]
This tower is asymptotically optimal.
\end{prop}

\proof
This is analogous to the proof of Proposition \ref{modulartowerprop}.
\Endproof

We illustrate in an example how easy it is to compute equations for
modular curves, by using the recursive definition.

\begin{exa}\label{X023exa}
We want to compute equations for the curve $X_0(2\cdot 3^m)$ in
characteristic $p\neq 2,3$. Our method is essentially the same as the
method of Elkies \cite{Elkies97}.

Note that the genus of $X_0(18)$ is zero. For $N=3,6,18$, we write
$L_N(u)=0$ for the differential equation corresponding to
$(H^1_{\dR}({\mathcal E}_N/S), \nabla)$.  Using the
description of the cusps in \cite{Shimura}, it is easy to check the
following statements. The differential equation $L_3(u)=0$ has three
singularities. It is no restriction to suppose that these
singularities are $x=0,1,\infty$, where $\infty$ maps to $j=0$ and
$x=0,1$ map to $j=\infty$ with ramification index $1,3,1$,
respectively. The map $X_0(6)\to X_0(3)$ of degree three is totally
branched above $x=\infty$ and has two points $P_0^1, P_0^2$ (resp.\
$P_1^1, P_1^2$ above $x=0$ (resp.\ $x=1$), where $P_\ast^e$ is
ramified of order $e$. Up to normalization, there is a unique such
cover which is given by $x=-27y^2/(y-4)^3$.  It follows that the
singularities of $L_6(u)=$ are $S_6=\{0,\infty,-8,1\}$. A look at the
ramification indices of these cusps in $X_0(6)\to X(1)$ tells us that
the Atkin--Lehner involution $\sigma_6$ acts on these points as
$(0,1)(-8, \infty)$, therefore $\sigma_6(y)=-8(x-1)/(x+8)$.

A similar argument shows that the natural projection $g:X_0(18)\to
X_0(6)$ is cyclic of order three and branched at $y=0, \infty$.
Therefore we may suppose that $X_0(18)\to X_0(6)$ is given by
$g(z)=z^3$. The singularities of $L_{18}(u)=0$ are just the
inverse image of $S_6$, i.e.\ $S_{18}=\{0, \infty, -2\zeta_3^i,
\zeta_3^i\}$, where $\zeta_3\in \FF_{p^2}$ is a primitive third
root of unity. The Atkin--Lehner involution is given, up to
normalization, by $\sigma_{18}(z)=-2(z-1)/(z+2)$. We define the
rational function $h(z)=\sigma_6\circ g\circ \sigma_{18}(z)=
z(z^2-2z+4)/(z^2+z+1)$.
%-(x-4)(7x^2-2x+4)/((2x+1)(x^2+x+7))$. KLOPT DIT WEL???
This gives the recursive definition for the modular curves
$X_0(2\cdot 3^m)$. \end{exa}

For $\lambda$ relatively prime to  $p$, we may define the
congruence subgroup $\Gamma_1(\lambda)\cap \Gamma_0(\ell^a)$. Write
$X_{1,0}(\lambda, \ell^a)_{\CC}$ for the quotient curve of the
completed upper half plane by $\Gamma_1(\lambda)\cap
\Gamma_0(\ell^a)$. It is well known that $X_{1,0}(\lambda,
\ell^a)_{\CC}$ has a model $X_{1,0}(\lambda, \ell^a)_{R}$ over
$R=\ZZ[\zeta_\lambda, 1/\lambda\ell]$, where $\zeta_\lambda$ is a
primitive $\lambda$th root of unity. Write
$\FF_q=\FF_p[\zeta_\lambda]$ and $X_{1,0}(\lambda,
\ell^a)=X_{1,0}(\lambda, \ell^a)_{R}\otimes \FF_q$. Let $f:X_{1,0}(\lambda,
\ell^a)\to X_0(\ell^a)$ be the natural projection.

  We can consider the pull back of
$(g,h):X_0(\ell^{a+1})\rightrightarrows X_0(\ell^a)$ via
$f:X_{1,0}(\lambda, \ell^a)\to X_0(\ell^a)$. Write ${\mathcal
T}_{g,h}(f):=(X_{1,0}(\lambda, \ell^{a+n})$ for the corresponding
tower. As in \cite{LMS}, one give a criterion on $p$ for this tower to
be asymptotically optimal. For example, it is easy to show that the
tower of Proposition \ref{exaprop} is isomorphic to $(X_{1,0}(8,
2^{4+m}))$. One can give an alternative proof of the facts on the
minimal splitting field by using this interpretation of the tower.

\flushright{ Insitut f\"ur Experimentelle Mathematik\\
Universit\"at Duisburg--Essen\\
Ellernstra{\ss}e 29\\ 45326 Essen\\
Germany\\
bouw@exp-math.uni-essen.de}

\vspace{-19ex}

\flushleft{
Fachbereich Mathematik\\
Universit\"at Duisburg--Essen\\
Universit\"atsstrasse 2\\
45117 Essen\\
Germany\\
peter.beelen@uni-essen.de}

\end{document}